\documentclass{amsart}

\newtheorem{theorem}{Theorem}
\newcommand{\bt}{\begin{theorem}}
\newcommand{\et}{\end{theorem}}

\newtheorem{lemma}{Lemma}
\newcommand{\bl}{\begin{lemma}}
\newcommand{\el}{\end{lemma}}

\newtheorem{corollary}{Corollary}
\newcommand{\bc}{\begin{corollary}}
\newcommand{\ec}{\end{corollary}}

\newcommand{\N}{\ensuremath{ \mathbf N }}

\newcommand{\mcm}{\ensuremath{ \mathcal M}}
\newcommand{\mcu}{\ensuremath{ \mathcal U}}

\newcommand{\beq}{\begin{equation}}
\newcommand{\eeq}{\end{equation}}
\newcommand{\benum}{\begin{enumerate}}
\newcommand{\eenum}{\end{enumerate}}

\title{Dirichlet series of integers with missing digits}
\author{Melvyn B. Nathanson}
\address{Lehman College (CUNY)}
\email{melvyn.nathanson@lehman.cuny.edu}



\date{\today}

\begin{document}
\maketitle

\begin{abstract}
For certain sequences $A$ of positive integers with missing $g$-adic digits, the Dirichlet series 
$F_A(s) = \sum_{a\in A} a^{-s}$ has abscissa of convergence $\sigma_c < 1$.  
The number $\sigma_c$ is computed.  
This generalizes and strengthens a classical theorem of Kempner 
on the convergence of the sum of the reciprocals 
of a sequence of integers with missing decimal digits.  
\end{abstract}

Let $A_{10,9}$ be the set of positive integers whose decimal representation contains no 9.  
A classical theorem of Kempner~\cite{kemp14} states that the harmonic series 
$\sum_{a\in A_{10,9}} 1/a$ converges.  
There is a straightforward generalization of this result.  Let $g$ be an integer such that $g \geq 2$.
Every integer $n$ in the interval $\left[ g^{m-1}, g^m -1 \right]$ has a unique $g$-adic representation 
\[
n = \sum_{i=0}^{m-1} c_i g^i 
\] 
with 
\[
c_i \in \{ 0,1,2,\ldots, g-1\} \qquad \text{for all $i \in \{0,1,\ldots, m-1\}$} 
\]
and 
\[
c_{m-1} \neq 0.
\]
Let $u \in [0,g-1]$ and let $A_{g,u}$ be the set of positive integers whose $g$-adic representation contains 
no digit $c_i = u$.  The series $ \sum_{a\in A} 1/a$ converges.  
This is Theorem 144 in Hardy and Wright~\cite{hard-wrig08}.  

It is natural to ask if Kempner's convergence theorem can be strengthened.  
Does there exist a real number $\sigma  < 1$ such that the infinite series $\sum_{a\in A_{g,u}} 1/a^{\sigma}$ converges?  A sharper question is: Compute the abscissa of convergence of the Dirichlet series 
\[
F_{A_{g,u}} (s) =  \sum_{a\in A_{g,u}} \frac{1}{a^s}. 
\]
We shall prove that this series has abscissa of convergence 
\[
\sigma_c = \frac{\log (g-1)}{\log g}
\] 
and that $F_{A_{g,u}}(\sigma_c)$ diverges.  
This is a corollary of Theorem~\ref{MissingDigits:theorem:2} below.

Let $\N = \{1,2,3,\ldots\}$ be the set of positive integers and 
$\N_0 = \{0,1,2,3,\ldots\}$ the set of nonnegative integers. 
For $x,y \in \N_0$, define the \emph{interval of integers}  $[x,y] = \{n\in \N_0: x \leq n \leq y\}$.

Fix an integer $g \geq 2$.  For all $i \in \N_0$, let $U_i$ be a proper subset of $[0,g-1]$, 
and let $\mcu = (U_i)_{i=0}^{\infty}$. 
Let $A_{g,\mcu}$ be the set of positive integers $n$ with $g$-adic representation $n = \sum_{i=0}^{m-1} c_ig^i$ 
such that $c_i \in [0,g-1] \setminus U_i$ for all $i \in [0,m-1]$ and $c_{m-1} \neq 0$.  
Consider the Dirichlet series 
\[
F_{A_{g,\mcu}}(s) =  \sum_{a\in A_{g,\mcu}} \frac{1}{a^s}. 
\]
This series converges if $\sigma = \Re(s) > 1$.
We have 
\[
A_{g,\mcu} \cap [g^{m-1}, g^m -1] \neq \emptyset
\]
if and only if $U_{m-1} \neq [1,g-1]$. 
Let 
\[
\mcm = \left\{m \in \N: U_{m-1} \neq [1,g-1] \right\}.
\]
The set $A_{g,\mcu}$ is infinite if and only if the set \mcm\ is infinite.  
If  \mcm\ is finite, then the series $F_A(s)$ is a Dirichlet polynomial, 
which is an entire function.  If \mcm\ is infinite, then $F_A(0)$ diverges, 
and so the Dirichlet series $F_{A_{g,\mcu}}(s)$ 
has abscissa of convergence $\sigma_c$ with $0 \leq \sigma_c \leq 1$.
We shall compute $\sigma_c$ for a large class of sets of integers with missing $g$-adic digits.

\bt          \label{MissingDigits:theorem:1}
For all $i \in \N_0$, let $U_i$ be a proper subset of $[0,g-1]$ such that 
(i) the set $\mcm = \left\{m \in \N: U_{m-1} \neq [1,g-1] \right\}$ is infinite, and 
(ii) 
there exist nonnegative real numbers $\alpha_0, \alpha_1, \ldots, \alpha_{g-1}$ 
such that, for all $m \in \N$ and $k \in [0,g-1]$, 
\beq                                        \label{MissingDigits:epsilon-alpha}
card\{i \in [0,m-1]: |U_i| = k\} = \alpha_k m+ \varepsilon_k(m)
\eeq
and   
\beq                                        \label{MissingDigits:epsilon-1}
 \lim_{m\rightarrow \infty} \frac{\varepsilon_k(m) }{m} = 0. 
\eeq
Let $\mcu = (U_i)_{i=0}^{\infty}$. 
Let $A_{g,\mcu}$  be the set of positive integers $n$ with $g$-adic representation 
$n = \sum_{i=0}^{m-1} c_ig^i$ 
such that  $c_{m-1} \neq 0$ and $c_i \notin U_i$ for all $i \in [0,m-1]$.
The Dirichlet series 
\[
F_{A_{g,\mcu}}(s)  = \sum_{a\in A_{g,\mcu}} \frac{1}{a^s}
\]
has abscissa of convergence 
\[
\sigma_c = \frac{1}{\log g}\sum_{k=0}^{g-1}\alpha_k  \log (g-k).
\]
\et

\begin{proof}
Let $A = A_{g,\mcu}$.     For all $m \in \N$, let  
\[
I_m = \left[ g^{m-1}, g^m - 1\right]   
\] 
and let  $n \in A \cap I_m$ have the $g$-adic representation $n = \sum_{i=0}^{m-1} c_ig^i$.  
For $i \in [0,m-2]$ there are $g-|U_{m-1}|$ choices for the digit $c_i$.
If $0 \in U_{m-1}$, there are $g-|U_{m-1}|$ choices for the digit $c_{m-1}$.
If $0 \notin U_{m-1}$, there are $g-|U_{m-1}|-1$ choices for the digit $c_{m-1}$.  
It follows that 
\beq                                        \label{MissingDigits:AIm-1}
|A \cap I_m| = 
\prod_{i=0}^{m-1} (g-|U_i|) \hspace{2.7cm} \text{if $0 \in U_{m-1}$} 
\eeq
and
\beq                                        \label{MissingDigits:AIm-2}
|A \cap I_m| = 
\left(\frac{g-|U_{m-1}|-1}{g-|U_{m-1}|}\right)   \prod_{i=0}^{m-1} (g-|U_i|) \qquad \text{if $0 \notin U_{m-1}$.} 
\eeq

Let 
 \[
\sigma >  \sigma_c = \frac{1}{\log g} \sum_{k=0}^{g-1} \alpha_k\log (g-k).
 \]
We shall prove that  the infinite series $F_A(\sigma)$ converges.  

Choose a real number $\delta$ such that 
\beq                                        \label{MissingDigits:delta-1}
0 <  \delta < \sigma -  \sigma_c  
=  \sigma - \frac{1}{\log g} \sum_{k=0}^{g-1} \alpha_k\log (g-k).
 \eeq
 Choose $m_0 = m_0(\delta)$ such that 
\beq                                        \label{MissingDigits:delta-m0}
| \varepsilon_k(m)| <  \left(\frac{\delta }{g}\right) m \qquad \text{for all $k \in [0,g-1]$ and $m \geq m_0$.} 
\eeq
Equations~\eqref{MissingDigits:epsilon-alpha},~\eqref{MissingDigits:AIm-1}, 
and~\eqref{MissingDigits:AIm-2} 
imply that, for $m \geq m_0$, we have 
\begin{align*}
|A\cap I_m| & \leq \prod_{i=0}^{m-1} (g-|U_i|) 
  = \prod_{k=0}^{g-1}  (g-k)^{\alpha_k m+ \varepsilon_k(m)} \\
&  \leq  g^{ \sum_{k=0}^{g-1} |\varepsilon_k(m)|}  \prod_{k=0}^{g-1}  (g-k)^{\alpha_k m} \\
& <   g^{  \delta m}  \prod_{k=0}^{g-1}  (g-k)^{\alpha_k m}  \\
& = \left(  g^{  \delta }  \prod_{k=0}^{g-1}  (g-k)^{\alpha_k } \right)^m.
\end{align*}
It follows that 
\begin{align*}
F_A(\sigma) & = \sum_{a\in A} \frac{1}{a^{\sigma}} 
=  \sum_{\substack{a\in A \\ a < g_{m_0-1}}} \frac{1}{a^{\sigma}} +  
\sum_{\substack{a\in A \\ a \geq g_{m_0-1}}} \frac{1}{a^{\sigma}} \\
&  = \sum_{\substack{a\in A \\ a < g_{m_0-1}}} \frac{1}{a^{\sigma}} 
 + \sum_{m= m_0}^{\infty} \sum_{a\in A\cap I_m} \frac{1}{a^{\sigma}}\\
& \leq  \sum_{\substack{a\in A \\ a < g_{m_0-1}}} \frac{1}{a^{\sigma}} + 
 \sum_{m= m_0}^{\infty} \frac{|A \cap I_m|}{g^{(m-1)\sigma}}\\
& \leq   \sum_{\substack{a\in A \\ a < g_{m_0-1}}} \frac{1}{a^{\sigma}}  
+ \sum_{m= m_0}^{\infty} \frac{ 
 \left(  g^{  \delta }  \prod_{k=0}^{g-1}  (g-k)^{\alpha_k } \right)^m
}{g^{(m-1)\sigma}} \\
& = \sum_{\substack{a\in A \\ a < g_{m_0-1}}} \frac{1}{a^{\sigma}}  
+  g^{\sigma}  \sum_{m= m_0}^{\infty} 
\left( \frac{ g^{  \delta }   \prod_{k=0}^{g-1} (g-k)^{\alpha_k } }{g^{\sigma}} \right)^m.  
\end{align*}
Inequality~\eqref{MissingDigits:delta-1} implies  
\[
0 <  \frac{ g^{  \delta }   \prod_{k=0}^{g-1} (g-k)^{\alpha_k } }{g^{\sigma}} < 1.
\]
and so the infinite series $F_A(\sigma)$ converges if $\sigma > \sigma_c$.

Let $\sigma < \sigma_c$.   We shall prove that the infinite series $F_A(\sigma)$ diverges.
Choose a real number $\delta$ such that 
\beq                                        \label{MissingDigits:delta-2}
0 <  \delta <  \sigma_c - \sigma = \frac{1}{\log g} \sum_{k=0}^{g-1} \alpha_k\log (g-k) - \sigma.
 \eeq
Let $m \geq m_0$.  
If $0 \leq |U_{m-1}| \leq g-2$, then 
\[
\frac{g-|U_{m-1}|-1}{g-|U_{m-1}|} \geq \frac{1}{2}.  
\]
From~\eqref{MissingDigits:epsilon-alpha},~\eqref{MissingDigits:AIm-1},~\eqref{MissingDigits:AIm-2},  
and~\eqref{MissingDigits:delta-m0}, we obtain 
\begin{align*}
|A\cap I_m| & \geq \left(\frac{g-|U_i|-1}{g-|U_i|}\right)   \prod_{i=0}^{m-1} (g-|U_i|)\\
& \geq  \frac{1}{2}   \prod_{i=0}^{m-1} (g-|U_i|)\\ 
&  =   \frac{1}{2}\prod_{k=0}^{g-1}  (g-k)^{\alpha_k m+ \varepsilon_k(m)} \\
& \geq  \frac{1}{2} g^{- \sum_{k=0}^{g-1} |\varepsilon_k(m)|}  \prod_{k=0}^{g-1}  (g-k)^{\alpha_k m} \\
& > \frac{1}{2} g^{ - \delta m}  \prod_{k=0}^{g-1}  (g-k)^{\alpha_k m}  \\
& =  \frac{1}{2}\left(  g^{ - \delta }  \prod_{k=0}^{g-1}  (g-k)^{\alpha_k } \right)^m.
\end{align*}
If $|U_{m-1}| = g-1$ and  $m \in \mcm$, then $0 \in U_{m-1}$  and $g-|U_{m-1}|=1$.  
It follows that 
\begin{align*}
|A\cap I_m| & = \prod_{i=0}^{m-1} (g-|U_i|) = \prod_{k=0}^{g-1}  (g-k)^{\alpha_k m+ \varepsilon_k(m)} \\
& \geq g^{- \sum_{k=0}^{g-1} |\varepsilon_k(m)|}  \prod_{k=0}^{g-1}  (g-k)^{\alpha_k m} \\
& >  g^{ - \delta m}  \prod_{k=0}^{g-1}  (g-k)^{\alpha_k m}  \\
& = \left(  g^{ - \delta }  \prod_{k=0}^{g-1}  (g-k)^{\alpha_k } \right)^m.
\end{align*}
If $n \in \N_0 \setminus \mcm$, then $U_{m-1} = [1,g-1]$ and $A\cap I_m = \emptyset$.  
We have 
\begin{align*}
F_A(\sigma) & = \sum_{a\in A} \frac{1}{a^{\sigma}} 
=  \sum_{\substack{a\in A \\ a < g_{m_1-1}}} \frac{1}{a^{\sigma}}  +   \sum_{\substack{m= m_1 \\ m \in \mcm}}^{\infty} \sum_{a\in A\cap I_m} \frac{1}{a^{\sigma}}\\
& \geq  
 \sum_{\substack{m= m_1 \\ m \in \mcm}}^{\infty} \frac{|A \cap I_m|}{g^{m \sigma}}
 >  \frac{1}{2} \sum_{\substack{m= m_1 \\ m \in \mcm}}^{\infty} 
 \frac{ \left(   g^{ - \delta }  \prod_{k=0}^{g-1}  (g-k)^{\alpha_k } \right)^m}{g^{m\sigma}} \\
& =  \frac{1}{2}  \sum_{\substack{m= m_1 \\ m \in \mcm}}^{\infty} 
\left( \frac{  g^{ - \delta }  \prod_{k=0}^{g-1} (g-k)^{\alpha_k } }{g^{\sigma}} \right)^m.
\end{align*}
Inequality~\eqref{MissingDigits:delta-2} implies  
\[
\frac{ g^{- \delta }   \prod_{k=0}^{g-1} (g-k)^{\alpha_k } }{g^{\sigma}} >  1
\]
and so the infinite series $F_A(\sigma)$ diverges if $\sigma <  \sigma_c$.    
This completes the proof.  
\end{proof}

\bc
Let $u_i \in [0,g-1]$ and $U_i = \{u_i\}$ for all $i \in \N_0$.  
Let $A=A_{g,\mcu}$ be the set of positive integers whose 
$g$-adic representation contains no digit $c_i = u_i$. 
The Dirichlet series $F_A(s) = \sum_{a\in A}a^{-s}$ has abscissa of convergence 
\[
\sigma_c = \frac{\log (g-1)}{\log g}.
\]
In particular, Kempner's series $F_{A_{10,9}}(s)$ has abscissa of convergence $\log 9/\log 10$.
\ec

\begin{proof}
Apply Theorem~\ref{MissingDigits:theorem:1} with $U_i = \{u_i\}$ for all $i \in \N_0$. 
We have $\alpha_0 = 0$, $\alpha_1 = 1$,  $\alpha_k = 0$ for all $k \in [2,g-1]$, and 
$\varepsilon_k(m) = 0$ for all $k \in [0,g-1]$.  

For Kempner's series, let $g=10$ and $u_i = 9$ for all $i \in \N_0$.  
\end{proof}

\bt          \label{MissingDigits:theorem:2}
For all $i \in \N_0$, let $U_i$ be a proper subset of $[0,g-1]$ such that (i)
the set $\mcm = \left\{m \in \N: U_{m-1} \neq [1,g-1] \right\}$ is infinite, and (ii)
there exist nonnegative real numbers $\alpha_0, \alpha_1, \ldots, \alpha_{g-1},$ and $\beta$ 
such that, for all $m \in \N$ and $k \in [0,g-1]$, 
\beq                                        \label{MissingDigits:epsilon-alpha-2}
card\{i \in [0,m-1]: |U_i| = k\} = \alpha_k m+ \varepsilon_k(m)
\eeq
and  
\beq                                        \label{MissingDigits:epsilon-2}
|\varepsilon_k(m)| < \beta. 
\eeq
Let $\mcu = (U_i)_{i=0}^{\infty}$. 
Let $A_{g,\mcu}$  be the set of positive integers $n$ with $g$-adic representation 
$n = \sum_{i=0}^{m-1} c_ig^i$ 
such that  $c_{m-1} \neq 0$ and $c_i \notin U_i$ for all $i \in [0,m-1]$.
The Dirichlet series 
\[
F_{A_{g,\mcu}}(s)  = \sum_{a\in A_{g,\mcu}} \frac{1}{a^s} 
\]
has abscissa of convergence 
\beq                                        \label{MissingDigits:sigma-c}
\sigma_c = \frac{1}{\log g}\sum_{k=0}^{g-1}\alpha_k  \log (g-k).
\eeq
Moreover, $F_{A_{g,\mcu}}(\sigma_c)$ diverges.  
\et

\begin{proof}
The only difference between Theorem~\ref{MissingDigits:theorem:1} 
and Theorem~\ref{MissingDigits:theorem:2}
is that condition~\eqref{MissingDigits:epsilon-1} has been replaced by 
the more restrictive condition~\eqref{MissingDigits:epsilon-2}.  
Thus, $F_{A_{g,\mcu}}(s)$ has abscissa of convergence $\sigma_c$.  
We shall prove that $F_{A_{g,\mcu}}(\sigma_c)$ diverges.  
Note that~\eqref{MissingDigits:sigma-c} is equivalent to 
\beq                                        \label{MissingDigits:sigma-c-1}
\frac{ \prod_{k=0}^{g-1} (g-k)^{\alpha_k } }{g^{ \sigma_c}} = 1
\eeq
and that~\eqref{MissingDigits:epsilon-2} implies 
\[
\sum_{k=0}^{g-1}| \varepsilon_k(m) | < g\beta
\]
for all $m \geq 1$.  
Following the proof of Theorem~\ref{MissingDigits:theorem:1}, 
we see that $|U_{m-1}| \leq g-2$ implies 
\begin{align*} |A_{g,\mcu} \cap I_m| & \geq  \frac{1}{2}   \prod_{i=0}^{m-1} (g-|U_i|)\\ 
&  =   \frac{1}{2}\prod_{k=0}^{g-1}  (g-k)^{\alpha_k m+ \varepsilon_k(m)} \\
& \geq  \frac{1}{2} g^{- \sum_{k=0}^{g-1} |\varepsilon_k(m)|}  \prod_{k=0}^{g-1}  (g-k)^{\alpha_k m} \\
& > \frac{1}{2g^{g\beta}}   \left( \prod_{k=0}^{g-1}  (g-k)^{\alpha_k } \right)^m.
\end{align*}
If $|U_{m-1}| = g-1$ and $0 \in U_{m-1}$, then  $g-|U_{m-1}| = 1$ and 
\begin{align*}
|A_{g,\mcu}  \cap I_m| & = \prod_{i=0}^{m-1} (g-|U_i|) = \prod_{k=0}^{g-1}  (g-k)^{\alpha_k m+ \varepsilon_k(m)} \\
& > \frac{1}{g^{g\beta}}  \left( \prod_{k=0}^{g-1}  (g-k)^{\alpha_k } \right)^m.
\end{align*}
Therefore, 
\begin{align*}
F_{A_{g,\mcu}}(\sigma_c) & = \sum_{a\in A_{g,\mcu} } \frac{1}{a^{\sigma_c}} 
\geq   \sum_{\substack{m= m_2 \\ m \in \mcm}}^{\infty} \sum_{a\in A_{g,\mcu}  \cap I_m} \frac{1}{a^{ \sigma_c}}\\
& \geq   \sum_{\substack{m= m_2 \\ m \in \mcm}}^{\infty} \frac{|A_{g,\mcu}  \cap I_m|}{g^{m  \sigma_c}} \\
&  > \frac{1}{2g^{g\beta}}  \sum_{\substack{m= m_2 \\ m \in \mcm}}^{\infty} 
\left( \frac{ \prod_{k=0}^{g-1} (g-k)^{\alpha_k } }{g^{ \sigma_c}} \right)^m.
\end{align*}
It follows from~\eqref{MissingDigits:sigma-c-1} that the infinite series 
$F_{A_{g,\mcu}}( \sigma_c)$ diverges.   
This completes the proof.  
\end{proof}

\bc
Let $u_i \in [0,g-1]$ and $U_i = \{u_i\}$ for all $i \in \N_0$.  
Let $A_{g,\mcu}$ be the set of positive integers whose 
$g$-adic representation contains no digit $c_i = u_i$. 
The Dirichlet series $F_{A_{g,\mcu}}(s) = \sum_{a\in A}a^{-s}$ has abscissa of convergence 
\[
\sigma_c = \frac{\log (g-1)}{\log g}
\]
and $F_{A_{g,\mcu} }(\sigma_c)$ diverges.  
In particular, Kempner's series $F_{A_{10,9}}(s)$ has abscissa of convergence $\log 9/\log 10$, 
and $F_{A_{10,9}}(\log 9/\log 10)$ diverges. 
\ec

\def\cprime{$'$} \def\cprime{$'$}
\providecommand{\bysame}{\leavevmode\hbox to3em{\hrulefill}\thinspace}
\providecommand{\MR}{\relax\ifhmode\unskip\space\fi MR }
\providecommand{\MRhref}[2]{%
  \href{http://www.ams.org/mathscinet-getitem?mr=#1}{#2}
}
\providecommand{\href}[2]{#2}

\end{document}